\documentclass[12pt, a4paper, oneside]{article}
\usepackage[top=1in, bottom=1.1in, left=0.9in, right=0.5in]{geometry}
\usepackage[parfill]{parskip}
\usepackage{amsmath}
\usepackage{amsthm}
\usepackage{graphicx}
\usepackage{amssymb}
\usepackage{epstopdf}
\usepackage{setspace}
\usepackage{authblk}
\usepackage{enumerate}
\usepackage{lmodern}
\usepackage[hypcap]{caption}
\usepackage[english]{babel}
\usepackage{fancyhdr}
\addto\captionsenglish{}
\addto\captionsenglish{}
\addto\captionsenglish{}
\addto\captionsenglish{}
\addto\captionsenglish{}

\newlength\longest

\newcommand{\cM}{{\mathcal M}}

\newcommand{\R}{{\mathbb R}}

\begin{document}

\newtheorem{theorem}{Theorem}[section]
\newtheorem{lemma}[theorem]{Lemma}
\newtheorem{corollary}[theorem]{Corollary}
\newtheorem{proposition}[theorem]{Proposition}
\newtheorem{conjecture}[theorem]{Conjecture}
\newtheorem{problem}[theorem]{Problem}
\newtheorem{claim}[theorem]{Claim}
\theoremstyle{definition}
\newtheorem{assumption}[theorem]{Assumption}
\newtheorem{remark}[theorem]{Remark}
\newtheorem{definition}[theorem]{Definition}
\newtheorem{example}[theorem]{Example}
\theoremstyle{remark}
\newtheorem{notation}{Notasi}
\renewcommand{\thenotation}{}

\title{On Inclusion Properties of Two Versions of Orlicz-Morrey Spaces}
\author{Al Azhary Masta${}^{1}$\footnote{\emph{Permanent Address}:
Department of Mathematics Education, Universitas Pendidikan Indonesia,
Jl. Dr. Setiabudi 229, Bandung 40154, Indonesia}, Hendra Gunawan${}^2$, Wono Setya-Budhi${}^3$}
\affil{${}^{1,2,3}$Analysis and Geometry Group, Faculty of Mathematics and\\
Natural Sciences, Bandung Institute of Technology,\\
Jl. Ganesha 10, Bandung 40132, Indonesia\\
E-mail: ${}^{1}$alazhari.masta@upi.edu, ${}^{2}$hgunawan@math.itb.ac.id,
${}^{3}$wono@math.itb.ac.id}
\date{}

\maketitle

\begin{abstract}
There are two versions of Orlicz-Morrey spaces (on $\mathbb{R}^n$), defined by Nakai in
2004 and by Sawano, Sugano, and Tanaka in 2012. In this paper we discuss the inclusion
properties of these two spaces and compare the results. Computing the norms of the
characteristic functions of balls in $\mathbb{R}^n$ is one of the keys to our results.
Similar results for weak Orlicz-Morrey spaces of both versions are also obtained.

\bigskip

\noindent{\bf Keywords}: Inclusion property, Orlicz-Morrey spaces, weak
Orlicz-Morrey spaces.\\
{\textbf{MSC 2010}}: Primary 46E30; Secondary 46B25, 42B35.
\end{abstract}

\section{Introduction}

Orlicz-Morrey spaces are generalizations of Orlicz spaces and Morrey spaces (on
$\mathbb{R}^n$).
There are two versions of Orlicz-Morrey spaces: one is defined by Nakai
\cite{Gala, Nakai1} and another by Sawano, Sugano, and Tanaka \cite{Gala, Sawano}.
We shall discuss both of them here. In particular, we are interested in the
inclusion properties of these spaces.

A function $\Phi:[0,\infty)\to [0,\infty)$ is called a Young function if $\Phi$
is convex, left-continuous, $\Phi(0) = 0$, and $ \lim \limits_{t\to\infty} \Phi(t)
= \infty$. Given two Young functions $\Phi, \Psi$, we write $\Phi \prec \Psi$ if
there exists a constant $C > 0$ such that $\Phi(t) \leq \Psi(Ct)$ for all $t > 0$.

Let $G_{1}$ be the set of all functions $\phi : (0, \infty) \rightarrow (0, \infty)$
such that $\phi(r)$ is nondecreasing but $\frac{\phi(r)}{r}$ is nonincreasing. For
a Young function $\Psi$, we also define $G_{2}$ to be the set of all functions $\psi :
(0, \infty) \rightarrow (0, \infty)$ such that $\psi(r)$ is nondecreasing but for
any $s > 0$, $\frac{\psi((r+s)^{n})}{\Psi^{-1}((\frac{r+s}{s})^{n})}$ is nonincreasing.

For $\phi_1,~\phi_2: (0,\infty) \rightarrow (0,\infty)$, we write $ \phi_1 \preceq \phi_2$
if there exists a constant $C > 0$ such that $\phi_1(t) \leq C \phi_2(t)$ for all $t > 0$.
If $\phi_1 \preceq \phi_2$ and $ \phi_2 \preceq \phi_1$, then we write $\phi_1 \approx
\phi_2$.

Let $\Phi$ be a Young function and $\phi \in G_{1}$. The Orlicz-Morrey spaces
$L_{\phi,\Phi}(\mathbb{R}^{n})$ (of Nakai's version) is the set of measurable functions
$f \in L^{1}_{\rm loc}(\mathbb{R}^{n})$ such that for every $a \in \mathbb{R}^n$ and $r > 0$,
the following quantity
$$\| f \|_{(\phi, \Phi, B(a,r))} := \inf \Bigl\{b>0 : \frac{\phi(|B(a,r)|)}
{|B(a,r)|} \int_{B(a,r)}\Phi\Bigl(\frac{|f(x)|}{b}\Bigr) dx \leq 1\Bigr\}$$
is finite. We use the notation $B(a,r)$ to denote the open ball in $\mathbb{R}^n$
centered at $a$ with radius $r$, and $|B(a,r)|$ for its Lebesgue measure. The
Orlicz-Morrey spaces $L_{\phi, \Phi}(\mathbb{R}^n)$ is a Banach space with respect
to the norm $\| f \|_{ L_{\phi,\Phi}(\mathbb{R}^{n})} :=
\sup\limits_{a \in \mathbb{R}^n,~r>0} \| f \|_{(\phi, \Phi, B(a,r))}.$

For $\phi(r)=r$, the space $L_{\phi,\Phi}(\mathbb{R}^{n})$ is the Orlicz space
$L_{\Phi}(\mathbb{R}^{n})$. Meanwhile, for $\Phi (r) = r^{p} $ and $\phi(r) =
r^{1 - \frac{\lambda}{n}}$ where $ 0 \leq \lambda \leq n$, the space
$L_{\phi,\Phi}(\mathbb{R}^{n})$ reduces to the Morrey space $L_{p,\lambda}(\mathbb{R}^{n})$.

Now, let $\Psi$ be a Young function and $\psi \in G_{2}$. Sawano, Sugano, and Tanaka
defined the Orlicz-Morrey space $\cM_{\psi,\Psi}(\mathbb{R}^{n})$ to be the set of measurable
functions $f \in L^{1}_{loc}(\mathbb{R}^{n})$ such that
$$\| f \|_{\cM_{\psi,\Psi}(\mathbb{R}^{n})} := \sup\limits_{a \in
\mathbb{R}^n,~r>0} \psi(|B(a,r)|) \| f \|_{(\Psi, B(a,r))} < \infty, $$
where $\| f \|_{(\Psi, B(a,r))} := \inf \bigl \{  {b>0:\frac{1}{|B(a,r)|}
\int_{B(a,r)}\Psi \bigl(\frac{|f(x)|}{b} \bigr) dx \leq1}\bigr \}$. Notice here that $\psi(|B(a,r)|)$ dominates of the growth $\|f\|_{(\Psi, B(a,r))}$. 

For the Young function $\Psi(t) = t^p$ ($1\leq p < \infty$), the spaces $\cM_{\psi,\Psi}(
\mathbb{R}^{n})$ are recognized as the generalized Morrey spaces $\cM^{p}_{\psi}(\mathbb{R}^{n})$.

Recently, Gunawan {\it et al.} \cite{Gunawan} presented a sufficient and necessary condition
for the inclusion relation between generalized Morrey spaces, as in the following theorem.

\bigskip

\begin{theorem}\label{thm:140812}
Let $1\le p_1\le p_2<\infty$ and $\psi_1, \psi_2 \in G_{2}$. Then the following statements are
equivalent:

\noindent{\rm (a)} $\psi_1 \preceq \psi_2$.

\noindent{\rm (b)} $\cM^{p_2}_{\psi_2}(\mathbb{R}^{n}) \subseteq \cM^{p_1}_{\psi_1}(\mathbb{R}^{n})$.

\noindent{\rm (c)} There exists a constant $C > 0$ such that $\|f\|_{\cM^{p_1}_{\psi_1}(
\mathbb{R}^{n})} \le C \|f\|_{\cM^{p_2}_{\psi_2}(\mathbb{R}^{n})}$ for every
$f\in \cM^{p_2}_{\psi_2}(\mathbb{R}^{n})$.
\end{theorem}

\medskip

In the same paper, Gunawan {\it et al.} also gave a necessary and sufficient condition for the
inclusion relation between generalized weak Morrey spaces.

Meanwhile, the inclusion relation between Orlicz spaces $L_{\Phi}(\R^n)$ and between weak
Orlicz spaces $wL_{\Phi}(\R^n)$ are known (see \cite{Lech, Masta1}).
In 2016, Masta {\it et al.} \cite{Masta2} also obtained the inclusion properties of Orlicz-Morrey
space $L_{\phi,\Phi}(\mathbb{R}^{n})$ of Nakai's version, as in the following theorem.

\bigskip

\begin{theorem}\label{theorem:1.7}
Let $\Phi_{1}, \Phi_{2}$ be Young functions and $\phi_{1}, \phi_{2} \in G_{1}$ such that
$\phi_{1} \approx \phi_{2}$. Then the following statements are equivalent:

{\parindent=0cm
{\rm (1)} $\Phi_{1} \prec \Phi_{2}$.

{\rm (2)} $L_{\phi_{2}, \Phi_{2}}(\mathbb{R}^n)  \subseteq L_{\phi_{1},\Phi_{1}}(\mathbb{R}^n)$.

{\rm (3)} There exists a constant $C>0$ such that $$\|f\|_{L_{\phi_{1},\Phi_{1}
}(\mathbb{R}^n)} \leq C \|f\|_{L_{\phi_{2}, \Phi_{2}}(\mathbb{R}^n)},$$ for every
$ f\in L_{\phi_{2}, \Phi_{2}}(\mathbb{R}^n) $.
\par}
\end{theorem}

\remark Note that the relation $\Phi_1 \prec \Phi_2$ is a necessary and sufficient condition
for the inclusion relation between Orlicz-Morrey spaces of Nakai's version.
For $\phi_{1}(t)=\phi_{2}(t)=t$, Theorem \ref{theorem:1.7} reduces to Theorem
3.4(a) in \cite{Lech}. Furthermore, for $\phi_{1}(t)=\phi_{2}(t)=t$ and $w_1 (x) = w_2 (x) = 1$,
Theorem \ref{theorem:1.7} complements Corollary 2.11 in \cite{Alen}, which states
that $\Phi_{1} \prec \Phi_{2}$ is a sufficient condition for inclusion relation between Orlicz spaces. Related results about inclusion properties of Orlicz-Morrey spaces can be found in \cite{Kita}.

In this paper, we would like to obtain the inclusion properties of
Orlicz-Morrey spaces $\cM_{\psi, \Psi}(\mathbb{R}^n)$ of Sawano-Sugano-Tanaka's version, and
compare it with the result for Nakai's version.
In addition, we will also prove similar results for weak Orlicz-Morrey spaces of both versions. With our results, we can see what parameters are significance in the inclusion properties for both versions.

To prove the results, we will use the same method as in \cite{Gunawan, Masta1, Masta2}, that is
by computing the norms of the characteristic functions of balls in $\R^n$. We also employ the
properties of the inverse function of $\Phi$, which are presented in the following lemma.

\bigskip

\begin{lemma}\label{lemma:1.1} \cite{Masta2,Nakai1,Rao}
Suppose that $\Phi$ is a Young function and $\Phi^{-1}$ denotes its inverse,
which is given by $\Phi^{-1}(s) :=\inf \{ r \geq 0 : \Phi (r) > s \}$ for
every $s \geq 0$. Then the followings hold:

{\parindent=0cm
{\rm (1)} $\Phi^{-1}(0) = 0$.

{\rm (2)} $ \Phi^{-1}(s_1) \leq \Phi^{-1}(s_2)$ for  $s_1 \leq s_2$.

{\rm (3)} $\Phi (\Phi^{-1}(s)) \leq s \leq \Phi^{-1}(\Phi(s))$ for $0 \leq s <
\infty$.

{\rm (4)} If, for some constants $C_1,C_2>0$, we have $\Phi_2^{-1}(s) \le C_1\Phi_1^{-1}(C_2s)$,
then $\Phi_1(\frac{t}{C_1}) \le C_2\Phi_2(t)$ for $t=\Phi_2^{-1}(s)$.

\par}
\end{lemma}

\medskip

Throughout this paper, the letter $C$ denotes a constant that may vary in values from
line to line. To keep track of some constants, we use subscripts, such as $C_{1}$ and
$C_{2}$.

\section{Inclusion Properties of Orlicz-Morrey Spaces}

As mentioned earlier, the key to our results is knowing the norms of the characteristic
balls in $\mathbb{R}^n$. Here is the first one on $\cM_{\psi, \Psi}(\mathbb{R}^n)$:

\bigskip

\begin{lemma}\label{lemma:1.11} \cite{Gala} For every $a\in\mathbb{R}^n$ and $r>0$, we have
$\| \chi_{B(a,r)} \|_{\cM_{\psi, \Psi}(\mathbb{R}^n)} =
\frac{\psi(|B(a,r)|)}{\Psi^{-1}(1)}$.
\end{lemma}

\bigskip

Our first theorem gives equivalent statements for the inclusion relation between
Orlicz-Morrey spaces of Sawano-Sugano-Tanaka's version.

\bigskip

\begin{theorem}\label{theorem:1.12}
Let $\Psi_1, \Psi_2$ be Young functions such that $\Psi_1 \prec \Psi_2$ and
$\psi_1, \psi_2 \in G_{2}$. Then the following statements are equivalent:

{\parindent=0cm
 {\rm (1)} $\psi_1 \preceq \psi_2$.

 {\rm (2)} $\cM_{\psi_2, \Psi_2}(\mathbb{R}^n) \subseteq \cM_{\psi_1, \Psi_1}(\mathbb{R}^n)$.

 {\rm (3)} There exists a constant $C > 0$ such that $$\| f \|_{\cM_{\psi_1,
 \Psi_1}(\mathbb{R}^n)} \leq C \| f \|_{\cM_{\psi_2, \Psi_2}(\mathbb{R}^n)}$$
 for every $ f \in \cM_{\psi_2, \Psi_2}(\mathbb{R}^n)$.

\par}
\end{theorem}

\noindent{\it Proof.} Let us first prove that (1) implies (2). Let $f \in \cM_{\psi_{2},
\Psi_{2}}(\mathbb{R}^{n})$. Recall that $\Psi_1 \prec \Psi_2$ means that there exists
a constant $C_1 > 0$ such that $\Psi_1(t) \leq \Psi_2(C_1t)$ for every $t > 0$.
For every $a\in\mathbb{R}^n$ and $r>0$, let $A_{(\Psi_1, B(a,r))} = \bigl \{  {b>0:\frac{1}{|B(a,r)|}
\int_{B(a,r)}\Psi_1 \bigl(\frac{|f(x)|}{C_1b} \bigr) dx \leq1}\bigr \}$
and $A_{(\Psi_2, B(a,r))} = \bigl \{  {b>0:\frac{1}{|B(a,r)|} \int_{B(a,r)}\Psi_2
\bigl(\frac{|f(x)|}{b} \bigr) dx \leq1}\bigr \}$.
Thus, for any $ b \in A_{(\Psi_2, B(a,r))}$, we have
\begin{align*}
 \frac{1}{|B(a,r)|} \int_{B(a,r)}\Psi_1 \left( \frac{|f(x)|}{C_1b} \right) dx &
 \leq \frac{1}{|B(a,r)|} \int_{B(a,r)} \Psi_2 \left(\frac{C_1|f(x)|}{C_1b} \right)dx\\
& = \frac{1}{|B(a,r)|} \int_{B(a,r)} \Psi_2 \left( \frac{|f(x)|}{b} \right) dx
\leq 1.
\end{align*}
Hence it follows that $ b \in A_{(\Psi_1, B(a,r))}$, and so we conclude that
$A_{(\Psi_2, B(a,r))} \subseteq A_{(\Psi_1, B(a,r))}$. Accordingly, we have
$$\left\|\frac{f}{C_1}\right\|_{(\Psi_1, B(a,r))} = \inf A_{(\Psi_1, B(a,r))}
\leq \inf A_{(\Psi_2, B(a,r))} = \| f \|_{(\Psi_2, B(a,r))},$$
and this holds for every $a\in\mathbb{R}^n$ and $r>0$.

Now there exists $C_2>0$ such that $\psi_1(s)\leq C_2 \psi_2(s)$ for every $s>0$. Combining
this with the previous estimate, we obtain
\begin{align*}
 \| f \|_{\cM_{\psi_1, \Psi_1} (\mathbb{R}^{n})} & =
 \sup\limits_{a \in \mathbb{R}^n,~r>0} \psi_1(|B(a,r)|) \| f \|_{(\Psi_1, B(a,r))}\\
 &\leq \sup\limits_{a \in \mathbb{R}^n,~r>0} C_1C_2 \psi_2(|B(a,r)|) \| f \|_{(\Psi_2, B(a,r))}\\
 & = C \| f \|_{\cM_{\psi_2, \Psi_2} (\mathbb{R}^{n})}.
\end{align*}
This proves that $\cM_{\psi_2, \Psi_2}(\mathbb{R}^{n}) \subseteq \cM_{\psi_1, \Psi_1} (\mathbb{R}^{n})$.

Next, since $(\cM_{\psi_2, \Psi_2}(\mathbb{R}^n), \cM_{\psi_1,\Psi_1}(\mathbb{R}^n))$ is a Banach pair,
it follows from \cite[Lemma 3.3]{Krein} that (2) and (3) are equivalent. It thus remains to show that
(3) implies (1).

Assume that (3) holds. Let $a\in\mathbb{R}^n$ and $ r > 0$. By Lemma \ref{lemma:1.11}, we have
$$ \frac{\psi_1(|B(a,r)|)}{\Psi^{-1}_{1}(1)} = \| \chi_{B(a,r)} \|_{\cM_{\psi_1, \Psi_1}(\mathbb{R}^n)}
\leq C \| \chi_{B(a,r)} \|_{\cM_{\psi_2, \Psi_2}(\mathbb{R}^n)} =  \frac{C \psi_2(|B(a,r)|)}{\Psi^{-1}_{2}(1)},$$
whence $\psi_1(|B(a,r)|) \leq \frac{C \Psi_{1}^{-1}(1)}{\Psi_{2}^{-1}(1)} \psi_2(|B(a,r)|)$.
Since $a\in\mathbb{R}^n$ and $r > 0$ are arbitrary, we get $\psi_1(t) \leq C_1 \psi_2(t)$ for every $t > 0$,
where $C_1 = \frac{C \Psi_{1}^{-1}(1)}{\Psi_{2}^{-1}(1)}$. \qed

\bigskip

\begin{corollary}\label{corollary:1.13}
Let $\Psi$ be a Young function and $\psi_1, \psi_2 \in G_{2}$. Then the following statements are equivalent:

{\parindent=0cm
 {\rm (1)} $\psi_1 \preceq \psi_2$.

 {\rm (2)} $\cM_{\psi_2, \Psi}(\mathbb{R}^n) \subseteq \cM_{\psi_1, \Psi}(\mathbb{R}^n)$.

 {\rm (3)} There exists a constant $C > 0$ such that $$\| f \|_{\cM_{\psi_1, \Psi}(\mathbb{R}^n)}
 \leq C \| f \|_{\cM_{\psi_2, \Psi}(\mathbb{R}^n)}$$ for every $f \in \cM_{\psi_2, \Psi}(\mathbb{R}^n)$.

\par}
\end{corollary}

\bigskip

\remark We note that the relation $\psi_1 \preceq \psi_2$ is a necessary and sufficient
condition for the inclusion relation between Orlicz-Morrey spaces of Sawano-Sugano-Tanaka's
version.

\section{Inclusion Properties of Weak Orlicz-Morrey Spaces}

We shall now discuss the inclusion properties of weak Orlicz-Morrey spaces. First, we recall
the definition of weak Orlicz-Morrey spaces $wL_{\phi, \Phi}(\mathbb{R}^n)$ \cite{Nakai2}.
Let $\Phi$ be a Young function and $\phi \in G_1$. The weak Orlicz-Morrey space
$wL_{\phi, \Phi}(\mathbb{R}^n)$ is the set of all measurable functions
$f : \mathbb{R}^n \rightarrow \mathbb{R} $ such that $\| f \|_{wL_{\phi, \Phi}(\mathbb{R}^n)}
= \sup\limits_{a \in \mathbb{R}^n, r > 0}\| f \|_{wL_{\phi, \Phi, B(a,r)}} < \infty$, where
\[
\| f \|_{wL_{\phi, \Phi, B(a,r)}} := \inf \left \{ b>0 :
\sup\limits_{t > 0} \frac{\Phi(t) \phi(|B(a,r)|) \bigl| \{ x \in B(a,r) : \frac{|f(x)|}{b} > t \}
\bigr|}{|B(a,r)|} \leq 1 \right \}
\]
for $a\in\mathbb{R}^n$ and $r>0$. If $\Phi(t) = t^p$, $1\leq p < \infty$ and $\phi(r) = r $, the space
$wL_{\phi, \Phi}(\mathbb{R}^n)$ is the weak Lebesgue space
$wL_{p}(\mathbb{R}^{n})$ (see \cite{Castillo}).

The relation between $wL_{\phi, \Phi}(\mathbb{R}^n)$ and $L_{\phi, \Phi}(\mathbb{R}^n)$ is
presented in the following lemma. (We leave the proof to the reader.)

\bigskip

\begin{lemma}\label{lemma:2.1}
Let $\Phi$ be a Young function and $\phi \in G_1$. Then $L_{\phi, \Phi}(\mathbb{R}^n)
\subseteq wL_{\phi, \Phi}(\mathbb{R}^n)$ with $\| f \|_{wL_{\phi, \Phi}(\mathbb{R}^n)}
\leq \| f \|_{L_{\phi, \Phi}(\mathbb{R}^n)}$ for every $f\in L_{\phi, \Phi}(\mathbb{R}^n)$.
\end{lemma}

\bigskip

The following lemma gives the norms of the characteristic functions of balls in $\mathbb{R}^n$.

\bigskip

\begin{lemma}\label{lemma:3.2}
Let $\Phi$ be a Young function, $\phi \in G_1$, $a\in\mathbb{R}^n$, and $r, r_{0} >0$. Then we have
\[\| \chi_{B(a,r_{0})} \|_{wL_{\phi, \Phi, B(a,r)}} = \frac{1}{\Phi^{-1}
\bigl( \frac{|B(a,r)|}{|B(a,r) \cap B(a,r_{0})| \phi(|B(a,r)|)} \bigr)}.\]
\end{lemma}

\noindent{\it Proof}.
Since $ \|\cdot\|_{wL_{\phi, \Phi, B(a,r)}} \leq \|\cdot\|_{(\phi, \Phi, B(a,r))}$
and  $\| \chi_{B(a,r_{0})} \|_{(\phi, \Phi, B(a,r))} = \frac{1}{\Phi^{-1}\bigl(
\frac{|B(a,r)|}{|B(a,r) \cap B(a,r_{0})|\phi(|B(a,r)|)} \bigr)}$, we obtain
$$\| \chi_{B(a,r_{0})} \|_{wL_{\phi, \Phi, B(a,r)}} \leq \frac{1}{\Phi^{-1} \bigl(
\frac{|B(a,r)|}{|B(a,r) \cap B(a,r_{0})|\phi(|B(a,r)|)} \bigr)}.$$
By the definitions of $\Phi^{-1}$ and $\| \cdot \|_{wL_{\phi, \Phi, B(a,r)}}$, we conclude that
$$\| \chi_{B(a,r_{0})} \|_{wL_{\phi, \Phi, B(a,r)}} = \frac{1}{\Phi^{-1} \bigl(
\frac{|B(a,r)|}{|B(a,r) \cap B(a,r_{0})| \phi(|B(a,r)|)} \bigr)}.$$ \qed

\begin{lemma}\label{lemma:3.3}
Let $\Phi$ be a Young function, $\phi \in G_1$, $a\in\mathbb{R}^n$, and $r_{0} >0$ be arbitrary. Then we have
$ \| \chi_{B(a,r_{0})} \|_{wL_{\phi, \Phi}(\mathbb{R}^n)} = \frac{1}{\Phi^{-1}(\frac{1}{ \phi(|B(a,r_{0})|)})}$.
\end{lemma}

\noindent{\it Proof}.
Since $\Phi$ is a Young function and $\phi \in G_1$, we have
$\| \chi_{B(a,r_{0})} \|_{L_{\phi, \Phi}(\mathbb{R}^n)} = \frac{1}{\Phi^{-1}(\frac{1}
{\phi(|B(a,r_{0})|)})}$ for $a\in\mathbb{R}^n$ and $r_0>0$ (see \cite{Gala}). Hence, by Lemma \ref{lemma:2.1},
we have $\| \chi_{B(a,r_{0})} \|_{wL_{\phi, \Phi}(\mathbb{R}^n)} \leq
\frac{1}{\Phi^{-1}(\frac{1}{ \phi(|B(a,r_{0})|)})}$. On the other hand,
\begin{align*}
\| \chi_{B(a,r_{0})} \|_{wL_{\phi, \Phi}(\mathbb{R}^n)}
& = \sup\limits_{a \in \mathbb{R}^n, r > 0}\| \chi_{B(a,r_{0})} \|_{wL_{\phi, \Phi, B(a,r)}}\\
& = \sup\limits_{a \in \mathbb{R}^n, r > 0} \frac{1}{\Phi^{-1} \bigl( \frac{|B(a,r)|}{|B(a,r)
\cap B(a,r_{0})| \phi(|B(a,r)|)} \bigr)}\\
& \geq \frac{1}{\Phi^{-1}(\frac{1}{ \phi(|B(a,r_{0})|)})}.
\end{align*}
Consequently, we have $ \| \chi_{B(a,r_{0})} \|_{wL_{\phi, \Phi}(\mathbb{R}^n)} =
\frac{1}{\Phi^{-1}(\frac{1}{ \phi(|B(a,r_{0})|)})}$. \qed

\medskip

Now we come to the inclusion property of weak Orlicz-Morrey spaces
of Nakai's version.

\bigskip

\begin{theorem}\label{theorem:2.4}
Let $\Phi_1, \Phi_2$ be Young functions, $\phi_1, \phi_2 \in G_{1}$ such that
$\phi_1 \preceq \phi_2$. Then the following statements are equivalent:

{\parindent=0cm
{\rm (1)} $\Phi_1 \prec \Phi_2$.

{\rm (2)} $wL_{\phi_2, \Phi_2}(\mathbb{R}^n) \subseteq wL_{\phi_1, \Phi_1}(\mathbb{R}^n)$.

{\rm(3)} There exists a constant $C > 0$ such that
$$\| f \|_{wL_{\phi_1, \Phi_1}(\mathbb{R}^n)} \leq C \| f \|_{wL_{\phi_2, \Phi_2} (\mathbb{R}^n)}$$
for every $ f \in wL_{\phi_2, \Phi_2}(\mathbb{R}^n)$.
\par}
\end{theorem}

\bigskip

\noindent{\it Proof}.
Assume that (1) holds and let $ f \in wL_{\phi_2, \Phi_2}(\mathbb{R}^n)$.
Since $\Phi_1 \prec \Phi_2$ and $\phi_1 \preceq \phi_2$, there exist constants
$C_1,C_2 > 0$ such that $\Phi_1(t) \leq \Phi_2(C_1t)$ and
$\phi_1(t)  \leq C_2 \phi_2(t)$ for every $ t > 0 $.
Let $a\in\mathbb{R}^n$ and $r>0$. We consider two cases.

Case I: $C_2 \geq 1$. Let
\begin{align*}
A_{\phi_1,\Phi_1, B(a,r)} &= \Bigl\{  {b>0: \sup\limits_{t > 0} \frac{\Phi_1(\frac{t}{C_2}) \phi_1(|B(a,r)|) \bigl|
\{ x \in B(a,r) : \frac{|f(x)|}{b} > t \} \bigr|}{|B(a,r)|} \leq1}\Bigr\} \\
& = \Bigl\{  {b>0: \sup\limits_{y > 0} \frac{\Phi_1(y) \phi_1(|B(a,r)|) \bigl|
\{ x \in B(a,r) : \frac{|f(x)|}{b} > C_2y \} \bigr|}{|B(a,r)|} \leq1}\Bigr\}\\
& = \Bigl\{  {b>0: \sup\limits_{y > 0} \frac{\Phi_1(y) \phi_1(|B(a,r)|) \bigl|
\{ x \in B(a,r) : \frac{|f(x)|}{C_2b} > y \} \bigr|}{|B(a,r)|} \leq1}\Bigr\}
\end{align*}
and
\begin{align*}
A_{\phi_2,\Phi_2, B(a,r)} & = \Bigl\{ b>0: \sup\limits_{t > 0} \frac{\Phi_2(C_1t) \phi_2(|B(a,r)|) \bigl|
\{ x \in B(a,r) : \frac{|f(x)|}{b} > t \} \bigr|}{|B(a,r)|} \leq 1 \Bigr\}\\
& = \Bigl\{ b>0: \sup\limits_{y > 0} \frac{\Phi_2(y) \phi_2(|B(a,r)|) \bigl|
\{ x \in B(a,r) : \frac{|f(x)|}{b} > \frac{y}{C_1} \} \bigr|}{|B(a,r)|} \leq 1\Bigr\}\\
& = \Bigl\{ b>0: \sup\limits_{y > 0} \frac{\Phi_2(y) \phi_2(|B(a,r)|) \bigl|
\{ x \in B(a,r) : \frac{|C_1 f(x)|}{b} > y \} \bigr|}{|B(a,r)|} \leq 1\Bigr\}.
\end{align*}
Then $\left\|\frac{f}{C_2}\right\|_{wL_{\phi_1, \Phi_1, B(a,r)}} = \inf A_{\phi_1, \Phi_1, B(a,r)}$
and $\| C_1f \|_{wL_{\phi_2, \Phi_2, B(a,r)}} = \inf A_{\phi_2, \Phi_2, B(a,r)}$.
Observe that, for arbitrary $ b \in A_{\phi_2, \Phi_2, B(a,r)}$ and $ t > 0$, we have
\begin{align*}
\frac{\Phi_1(\frac{t}{C_2}) \phi_1(|B(a,r)|) \bigl| \{ x \in B(a,r) : \frac{|f(x)|}{b} > t \}
\bigr|}{|B(a,r)|} & \leq \frac{\Phi_1(t) \frac{\phi_1(|B(a,r)|)}{C_2} \bigl|
\{ x \in B(a,r) : \frac{|f(x)|}{b} > t \} \bigr|}{|B(a,r)|}\\
& \leq \frac{\Phi_2(C_1t) \phi_2(|B(a,r)|) \bigl| \{ x \in B(a,r) : \frac{|f(x)|}{b} > t \} \bigr|}{|B(a,r)|}\\
& \leq \sup\limits_{y > 0} \frac{\Phi_2(y) \phi_2(|B(a,r)|) \bigl| \{ x \in B(a,r) : \frac{|C_1f(x)|}{b} > y \} \bigr|}{|B(a,r)|}\\
& \leq 1.
\end{align*}
Since $ t > 0$ is arbitrary, we have $\sup\limits_{t > 0} \frac{\Phi_1(\frac{t}{C_2}) \phi_1(|B(a,r)|)
\bigl| \{ x \in B(a,r) : \frac{|f(x)|}{b} > t \} \bigr|}{|B(a,r)|} \leq 1$.
Hence it follows that $b \in A_{\phi_1, \Phi_1, B(a,r)}$, and so we conclude that
$A_{\phi_2,\Phi_2, B(a,r)} \subseteq A_{\phi_1, \Phi_1, B(a,r)}$. Accordingly, we obtain
$$\left\|\frac{f}{C_2}\right\|_{wL_{\phi_1,\Phi_1, B(a,r)}} = \inf A_{\phi_1, \Phi_1, B(a,r)}
\leq \inf A_{\phi_2, \Phi_2, B(a,r)} = \| C_1f \|_{wL_{\phi_2, \Phi_2, B(a,r)}}.$$

Case II: $0 < C_2 < 1$. Observe that, for arbitrary $ b \in A_{(\phi_2, \Phi_2, B(a,r))}$
and $ t > 0$, we have
\begin{align*}
\frac{\Phi_1(t) \phi_1(|B(a,r)|) \bigl| \{ x \in B(a,r) : \frac{|C_2 f(x)|}{b} > t \}
\bigr|}{|B(a,r)|} & = \frac{\Phi_1(t) \phi_1(|B(a,r)|) \bigl| \{ x \in B(a,r) :
\frac{|f(x)|}{b} > \frac{t}{C_2} \} \bigr|}{|B(a,r)|} \\
& = \frac{\Phi_1(C_2t_1) \phi_1(|B(a,r)|) \bigl| \{ x \in B(a,r) : \frac{|f(x)|}{b} > t_1 \} \bigr|}{|B(a,r)|} \\
& \leq \frac{C_2 \Phi_1(t_1) \phi_1(|B(a,r)|) \bigl| \{ x \in B(a,r) : \frac{|f(x)|}{b} > t_1 \} \bigr|}{|B(a,r)|}\\
& \leq \frac{C_2^2 \Phi_2(C_1t_1) \phi_2(|B(a,r)|) \bigl| \{ x \in B(a,r) : \frac{|f(x)|}{b} > t_1 \} \bigr|}{|B(a,r)|}\\
& \leq \frac{\Phi_2(C_1t_1) \phi_2(|B(a,r)|) \bigl| \{ x \in B(a,r) : \frac{|f(x)|}{b} > t_1 \} \bigr|}{|B(a,r)|}\\
& \leq \sup\limits_{y > 0} \frac{\Phi_2(y) \phi_2(|B(a,r)|) \bigl| \{ x \in B(a,r) : \frac{|C_1f(x)|}{b} > y \} \bigr|}{|B(a,r)|}\\
& \leq 1.
\end{align*}
for $t_1 = \frac{t}{C_2}$. Since $ t > 0$ is arbitrary, we have $\sup\limits_{t > 0} \frac{\Phi_1(t) \phi_1(|B(a,r)|)
\bigl| \{ x \in B(a,r) : \frac{|C_2f(x)|}{b} > t \} \bigr|}{|B(a,r)|} \leq 1$.
Accordingly, we obtain $\| C_2 f\|_{wL_{\phi_1,\Phi_1, B(a,r)}} \leq  \| C_1f \|_{wL_{\phi_2, \Phi_2, B(a,r)}}.$

From Cases I and II, there exists a constant $C > 0$ such that $\| f\|_{wL_{\phi_1,\Phi_1, B(a,r)}}
\leq C \| f \|_{wL_{\phi_2, \Phi_2, B(a,r)}}.$
Since $a \in \mathbb{R}^n$ and $ r > 0$ are arbitrary, we conclude that
$\| f\|_{wL_{\phi_1,\Phi_1}(\mathbb{R}^n)} \leq C \| f\|_{wL_{\phi_2,\Phi_2}(\mathbb{R}^n)}$,
which implies that $wL_{\phi_2, \Phi_2}(\mathbb{R}^n) \subseteq wL_{\phi_1, \Phi_1}(\mathbb{R}^n)$.

As mentioned in \cite[Appendix G]{Nakai2}, we know that Lemma 3.3 in \cite{Krein} still holds
for quasi-Banach spaces, so (2) and (3) are equivalent.

Now, we will show that (3) implies (1). To do so, assume that (3) holds. By Lemma \ref{lemma:3.3}, we have
$$
 \frac{1}{\Phi^{-1}_{1}(\frac{1}{\phi_1(|B(a,r_{0})|)})} = \| \chi_{B(a,r_{0})} \|_{wL_{\phi_1, \Phi_1}(\mathbb{R}^n)}
 \leq C \| \chi_{B(a,r_{0})} \|_{wL_{\phi_2, \Phi_2}(\mathbb{R}^n)} =  \frac{C}{\Phi^{-1}_{2}(\frac{1}{\phi_2(|B(a,r_{0})|)})},
$$
whence $\Phi^{-1}_{2}(\frac{1}{\phi_1(|B(a,r_{0})|)}) \le C\Phi^{-1}_{1}(\frac{1}{\phi_{2}(|B(a,r_{0})|)})
\le C\Phi_1^{-1}(\frac{C_2}{\phi_1(|B(a,r_{0})|)})$,
for arbitrary $a\in\mathbb{R}^n$ and $r_0 > 0$. By Lemma \ref{lemma:1.1}(4), we have
$$
\Phi_1\Bigl(\frac{t_0}{C}\Bigr) \le C_2\Phi_2(t_0),
$$
where $t_0=\Phi_2^{-1}(\frac{1}{\phi_2(|B(a,r_{0})|)})$. If $C_2\le 1$, then
$
\Phi_1(\frac{t_0}{C}) \le \Phi_2(t_0).
$
If $C_2>1$, then noting that $\Phi_1$ is convex, we have
\[
\Phi_1\Bigl(\frac{t_0}{C_2C}\Bigr)\le \frac{1}{C_2}\Phi_1\Bigl(\frac{t_0}{C}\Bigr)\le \Phi_2(t_0).
\]
Since $a\in\mathbb{R}^n$ and $r_{0} > 0$ are arbitrary, we conclude that there exists $C_3>0$
such that
$\Phi_1(\frac{t}{C_3})\le \Phi_2(t)$ or equivalently $\Phi_1(t) \leq \Phi_2(C_3t)$
for every $t > 0$. \qed

\bigskip

\remark For $\phi_1(t)=\phi_2(t)= t$, Theorem \ref{theorem:2.4} reduces to Theorem 3.3 in \cite{Masta1}.

\bigskip

We shall now study the inclusion properties of weak Orlicz-Morrey spaces $w\cM_{\psi, \Psi}(\mathbb{R}^n)$.
Let $\Psi$ be a Young function and $\psi \in G_{2}$. The weak Orlicz-Morrey space
$w\cM_{\psi, \Psi}(\mathbb{R}^n)$ is the set of all measurable functions $f : \mathbb{R}^n \rightarrow
\mathbb{R}$ such that
$$
\| f \|_{w\cM_{\psi, \Psi}(\mathbb{R}^n)}  = \sup\limits_{a \in \mathbb{R}^n, r > 0} \psi(|B(a,r)|)
\| f \|_{w\cM_{\Psi, B(a,r)}} < \infty
$$
where $$\| f \|_{w\cM_{\Psi, B(a,r)}} := \inf \left \{ b>0:\sup\limits_{t > 0} \frac{\Psi(t) \Bigl|
\{ x \in B(a,r) : \frac{|f(x)|}{b} > t \}\Bigr|}{|B(a,r)|} \leq 1\right \},$$
for $a\in\mathbb{R}^n$ and $r>0$. Note that if there exists $C>0$ such that $\Psi_1(t)\le \Psi_2(Ct)$
for every $t>0$, then $\|f\|_{w\cM_{\Psi_1, B(a,r)}} \le C\|f\|_{w\cM_{\Psi_2, B(a,r)}}$ for every
$a\in\mathbb{R}^n$ and $r>0$.

The following lemma tells us that $w\cM_{\psi, \Psi}(\mathbb{R}^n)$ contains $\cM_{\psi, \Psi}(\mathbb{R}^n)$.
(We leave the proof to the reader.)

\bigskip

\begin{lemma}\label{lemma:3.6}
Let $\Psi$ be a Young function and $\psi \in G_2$. Then $\cM_{\psi, \Psi}(\mathbb{R}^n) \subseteq
w\cM_{\psi, \Psi}(\mathbb{R}^n)$ with $\| f \|_{w\cM_{\psi, \Psi}(\mathbb{R}^n)} \leq
\| f \|_{\cM_{\psi, \Psi}(\mathbb{R}^n)}$ for every $f\in \cM_{\psi, \Psi}(\mathbb{R}^n)$.
\end{lemma}

\bigskip

Following similar arguments as in the proof of Lemma \ref{lemma:3.2}, we have the following lemma.

\bigskip

\begin{lemma}\label{lemma:3.7}
Let $\Psi$ be a Young function, $a\in\mathbb{R}^n$, and $r, r_{0} >0$. Then we have
$$ \| \chi_{B(a,r_{0})} \|_{w\cM_{\Psi, B(a,r)}} = \frac{1}{\Psi^{-1} \bigl(
\frac{|B(a,r)|}{|B(a,r) \cap B(a,r_{0})|} \bigr)}.$$
\end{lemma}

\bigskip

The norms of the characteristic functions of balls in $\mathbb{R}^n$ is presented in the following
lemma.

\bigskip

\begin{lemma}\label{lemma:3.8}
Let $\Psi$ be a Young function, $\psi\in G_2$, $a\in\mathbb{R}^n$, and $r_{0} >0$ be arbitrary. Then we have
$ \| \chi_{B(a,r_{0})} \|_{w\cM_{\psi, \Psi}(\mathbb{R}^n)} = \frac{\psi(|B(a,r_{0})|)}{\Psi^{-1}(1)}$.
\end{lemma}

\noindent{\it Proof}.
Since $\Psi$ is a Young function and $\psi \in G_2$, by Lemma \ref{lemma:1.11} and Lemma \ref{lemma:3.6}
we have $ \| \chi_{B(a,r_{0})} \|_{w\cM_{\psi, \Psi}(\mathbb{R}^n)} \leq
\frac{\psi(|B(a,r_{0})|)}{\Psi^{-1}(1)}$, for $a\in\mathbb{R}^n$ and $r_0>0$. On the other hand, we have
\begin{align*}
\| \chi_{B(a,r_{0})} \|_{w\cM_{\psi, \Psi}(\mathbb{R}^n)}  & = \sup\limits_{a \in \mathbb{R}^n, r > 0}
\psi(|B(a,r)|)\| \chi_{B(a,r_{0})} \|_{w\cM_{\Psi, B(a,r)}}\\
& =  \sup\limits_{a \in \mathbb{R}^n, r > 0} \frac{\psi(|B(a,r)|)}{\Psi^{-1} \bigl(
\frac{|B(a,r)|}{|B(a,r) \cap B(a,r_{0})|} \bigr)}
 \geq \frac{\psi(|B(a,r_0)|)}{\Psi^{-1}(1)}.
\end{align*}
Consequently, we have $ \| \chi_{B(a,r_{0})} \|_{w\cM_{\psi, \Psi}(\mathbb{R}^n)} =
\frac{\psi(|B(a,r_{0})|)}{\Psi^{-1}(1)}$, as desired. \qed

\bigskip

Now we come to the inclusion property of weak Orlicz-Morrey spaces of Sawano-Sugano-Tanaka's version.

\bigskip

\begin{theorem}\label{theorem:2.9}
Let $\Psi_1, \Psi_2$ be Young functions such that $\Psi_1 \prec \Psi_2 $ and
$\psi_1, \psi_2 \in G_{2}$. Then the following statements are equivalent:

{\parindent=0cm
{\rm (1)} $\psi_1 \preceq \psi_2$.

{\rm (2)} $w\cM_{\psi_2, \Psi_2}(\mathbb{R}^n) \subseteq w\cM_{\psi_1, \Psi_1}(\mathbb{R}^n)$.

{\rm (3)} There exists a constant $C > 0$ such that $$\| f \|_{w\cM_{\psi_1, \Psi_1}(\mathbb{R}^n)}
\leq \| f \|_{w\cM_{\psi_2, \Psi_2} (\mathbb{R}^n)}$$ for every $ f \in w\cM_{\psi_2, \Psi_2}(\mathbb{R}^n)$.
\par}
\end{theorem}

\noindent{\it Proof}.
Assume that (1) holds. Let $ f \in w\cM_{\psi_2,\Psi_2}(\mathbb{R}^n)$. Since $\Psi_1 \prec \Psi_2$
and $\psi_1 \preceq \psi_2$,  there exist constant $C_1, C_2 > 0$ such that $\psi_1(t) \leq C_1\psi_2(t)$
and $\Psi_1(t) \leq \Psi_2(C_2t)$ for every $t > 0$. Observe that
\begin{align*}
 \| f \|_{w\cM_{\psi_1, \Psi_1}(\mathbb{R}^n)} & = \sup\limits_{a \in \mathbb{R}^n, r>0}
 \psi_1(|B(a,r)|) \| f \|_{w\cM_{\Psi_1, B(a,r)}}\\
 &\leq \sup\limits_{a \in \mathbb{R}^n, r>0} C_1 \psi_2(|B(a,r)|) \| f \|_{w\cM_{\Psi_1, B(a,r)}}\\
 &\leq \sup\limits_{a \in \mathbb{R}^n, r>0} C_1C_2\psi_2(|B(a,r)|) \| f \|_{w\cM_{\Psi_2, B(a,r)}}\\
 & = C_1C_2 \| f \|_{w\cM_{\psi_2, \Psi_2}(\mathbb{R}^n)}.
\end{align*}
Hence we conclude that $w\cM_{\psi_2, \Psi_2}(\mathbb{R}^n) \subseteq w\cM_{\psi_1, \Psi_1}(\mathbb{R}^n)$.

As before, (2) and (3) are equivalent, and so it remains to show that (3) implies (1). To do so, assume that
(3) holds. By Lemma \ref{lemma:3.8}, we have
$$
 \frac{\psi_1(|B(a,r)|)}{\Psi^{-1}_{1}(1)} = \| \chi_{B(a,r)} \|_{w\cM_{\psi_1, \Psi_1}(\mathbb{R}^n)}
 \leq C  \| \chi_{B(a,r)} \|_{w\cM_{\psi_2, \Psi_2}(\mathbb{R}^n)} =  \frac{C \psi_2 (|B(a,r)|)}{\Psi^{-1}_{2}(1)},
$$
whence $ \psi_1 (|B(a,r)|) \leq \frac{C\Psi^{-1}_{1}(1)}{\Psi^{-1}_{2}(1)} \psi_2 (|B(a,r)|)$,
for every $a\in\mathbb{R}^n$ and $r > 0$. We conclude that
$$\psi_1(t) \leq C_1 \psi_2(t)$$ for every $t > 0$, with $C_1 = \frac{C\Psi^{-1}_{1}(1)}{\Psi^{-1}_{2}(1)}$. \qed

\bigskip


\section{Further Results and Concluding Remarks}

The inclusion properties of Orlicz-Morrey spaces $L_{\phi,\Phi}(\mathbb{R}^n)$ (Theorem \ref{theorem:1.7})
and weak Orlicz-Morrey spaces $wL_{\phi,\Phi}(\mathbb{R}^n)$ (Theorem \ref{theorem:2.4}) generalize the
inclusion properties of Orlicz spaces and weak Orlicz spaces in \cite{Lech, Masta1}. Meanwhile, the inclusion
properties of Orlicz-Morrey spaces $\cM_{\psi,\Psi}(\mathbb{R}^n)$ (Theorem \ref{theorem:1.12}) and weak
Orlicz-Morrey spaces $w\cM_{\psi,\Psi}(\mathbb{R}^n)$ (Theorem \ref{theorem:2.9}) generalize the
inclusion properties of generalized Morrey spaces and generalized weak Morrey spaces in \cite{Gunawan}. Combining the results, one realizes that the inclusion relation between Orlicz-Morrey spaces is equivalent to that between weak Orlicz-Morrey spaces.

\medskip

Recently, Guliyev,\textit{ et al.} \cite{Deringoz,Guliyev} also introduced (strong) Orlicz-Morrey spaces different from Nakai's or Sawano-Sugano-Tanaka's versions. For a Young function $\Theta$, let $G_{\Theta}$ be the set of all functions $\theta : (0, \infty) \rightarrow (0, \infty)$
such that $\theta(r)$ is decreasing but $\Theta^{-1}(t^{-n})\theta(t)^{-1}$ is almost decreasing for all $ t>0$. Now, let $\Theta$ be a Young function and $\theta \in G_{\Theta}$. As in \cite{Guliyev}, may they define the Orlicz-Morrey space $\cM_{\theta,\Theta}(\mathbb{R}^{n})$ to be the set of measurable
functions $f$ such that
$$\| f \|_{\cM_{\theta,\Theta}(\mathbb{R}^{n})} := \sup\limits_{a \in
	\mathbb{R}^n,~r>0} \frac{1}{\theta(|B(a,r)|^{\frac{1}{n}})}\Theta^{-1}\Bigl(\frac{1}{|B(a,r)|}\Bigr) \|f\|_{L_{\Theta}(B(a,r))} < \infty, $$
where $\| f \|_{L_{\Theta}(B(a,r))} := \inf \bigl \{  {b>0:
	\int_{B(a,r)}\Theta \bigl(\frac{|f(x)|}{b} \bigr) dx \leq1}\bigr \}$.

By using similar arguments in the proof of Theorem \ref{theorem:1.7} and Theorem \ref{theorem:1.12}, one may obtain the following theorem.

\medskip

\begin{theorem}\label{theorem:4.2}
Let $\Theta_1, \Theta_2$ be Young functions such that $\Theta_1 \prec \Theta_2$, $\Theta^{-1}_1 \prec \Theta^{-1}_2$, and $\theta_1, \theta_2 \in G_{\Theta}$. Then the following statements are equivalent:

{\parindent=0cm
	{\rm (1)} $\theta_2 \preceq \theta_1$.
	
	{\rm (2)} $\cM_{\theta_2, \Theta_2}(\mathbb{R}^n) \subseteq \cM_{\theta_1, \Theta_1}(\mathbb{R}^n)$.
	
	{\rm(3)} There exists a constant $C > 0$ such that
	$$\| f \|_{\cM_{\theta_1, \Theta_1}(\mathbb{R}^n)} \leq C \| f \|_{\cM_{\theta_2, \Theta_2} (\mathbb{R}^n)}$$
	for every $ f \in \cM_{\theta_2, \Theta_2}(\mathbb{R}^n)$.
	\par}
\end{theorem}

Comparing Theorem \ref{theorem:1.12} and Theorem \ref{theorem:4.2}, we can say that the condition on the growth parameters for the inclusion of Orlicz-Morrey spaces $\cM_{\psi,\Psi}(\mathbb{R}^n)$ and $\cM_{\theta,\Theta}(\mathbb{R}^n)$ are in principal the same. However, the condition on the Young function for the inclusion of the Orlicz-Morrey space $\cM_{\psi,\Psi}(\mathbb{R}^n)$ is simpler than that for the Orlicz-Morrey space  $\cM_{\theta,\Theta}(\mathbb{R}^n)$.

\medskip

\textbf{Acknowledgement}. The first and second authors are supported by ITB Research and Innovation 
Program No. 006d/I1.C01/PL/2016.

\end{document}